\documentclass[11pt]{article}

\begin{document}



\newcommand{\newc}{\newcommand}


\renewcommand{\theequation}{\thesection.\arabic{equation}}
\newc{\eqnoset}{\setcounter{equation}{0}}
\newcommand{\myref}[2]{#1~\ref{#2}}

\newcommand{\mref}[1]{(\ref{#1})}
\newcommand{\reflemm}[1]{Lemma~\ref{#1}}
\newcommand{\refrem}[1]{Remark~\ref{#1}}
\newcommand{\reftheo}[1]{Theorem~\ref{#1}}
\newcommand{\refdef}[1]{Definition~\ref{#1}}
\newcommand{\refcoro}[1]{Corollary~\ref{#1}}
\newcommand{\refprop}[1]{Proposition~\ref{#1}}
\newcommand{\refsec}[1]{Section~\ref{#1}}
\newcommand{\refchap}[1]{Chapter~\ref{#1}}

\newcommand{\beq}{\begin{equation}}
\newcommand{\eeq}{\end{equation}}
\newcommand{\beqno}[1]{\begin{equation}\label{#1}}

\newcommand{\barr}{\begin{array}}
\newcommand{\earr}{\end{array}}

\newc{\bearr}{\begin{eqnarray*}}
\newc{\eearr}{\end{eqnarray*}}

\newc{\bearrno}[1]{\begin{eqnarray}\label{#1}}
\newc{\eearrno}{\end{eqnarray}}

\newc{\non}{\nonumber}
\newc{\nol}{\nonumber\nl}

\newcommand{\bdes}{\begin{description}}
\newcommand{\edes}{\end{description}}
\newc{\benu}{\begin{enumerate}}
\newc{\eenu}{\end{enumerate}}
\newc{\btab}{\begin{tabular}}
\newc{\etab}{\end{tabular}}



\newtheorem{theorem}{Theorem}[section]
\newtheorem{defi}[theorem]{Definition}
\newtheorem{lemma}[theorem]{Lemma}
\newtheorem{rem}[theorem]{Remark}
\newtheorem{exam}[theorem]{Example}
\newtheorem{propo}[theorem]{Proposition}
\newtheorem{corol}[theorem]{Corollary}

\renewcommand{\thelemma}{\thesection.\arabic{lemma}}

\newcommand{\btheo}[1]{\begin{theorem}\label{#1}}
\newc{\brem}[1]{\begin{rem}\label{#1}\em}
\newc{\bexam}[1]{\begin{exam}\label{#1}\em}
\newc{\bdefi}[1]{\begin{defi}\label{#1}}
\newcommand{\blemm}[1]{\begin{lemma}\label{#1}}
\newcommand{\bprop}[1]{\begin{propo}\label{#1}}
\newcommand{\bcoro}[1]{\begin{corol}\label{#1}}
\newcommand{\etheo}{\end{theorem}}
\newcommand{\elemm}{\end{lemma}}
\newcommand{\eprop}{\end{propo}}
\newcommand{\ecoro}{\end{corol}}
\newc{\erem}{\end{rem}}
\newc{\eexam}{\end{exam}}
\newc{\edefi}{\end{defi}}

\newc{\rmk}[1]{{\bf REMARK #1: }}
\newc{\DN}[1]{{\bf DEFINITION #1: }}

\newcommand{\bproof}{{\bf Proof:~~}}
\newc{\eproof}{{\vrule height8pt width5pt depth0pt}\vspace{3mm}}


\newcommand{\rarr}{\rightarrow}
\newcommand{\Rarr}{\Rightarrow}
\newcommand{\tru}{\backslash}
\newc{\bfrac}[2]{\dspl{\frac{#1}{#2}}}


\newc{\nl}{\vspace{2mm}\\}
\newc{\nid}{\noindent}


\newcommand{\oneon}[1]{\frac{1}{#1}}
\newcommand{\dspl}{\displaystyle}
\newc{\grad}{\nabla}
\newc{\Div}{\mbox{div}}
\newc{\pdt}[1]{\dspl{\frac{\partial{#1}}{\partial t}}}
\newc{\pdn}[1]{\dspl{\frac{\partial{#1}}{\partial \nu}}}
\newc{\pdNi}[1]{\dspl{\frac{\partial{#1}}{\partial \mathcal{N}_i}}}
\newc{\pD}[2]{\dspl{\frac{\partial{#1}}{\partial #2}}}
\newc{\dt}{\dspl{\frac{d}{dt}}}
\newc{\bdry}[1]{\mbox{$\partial #1$}}
\newc{\sgn}{\mbox{sign}}

\newc{\Hess}[1]{\frac{\partial^2 #1}{\pdh z_i \pdh z_j}}
\newc{\hess}[1]{\partial^2 #1/\pdh z_i \pdh z_j}


\newcommand{\Coone}[1]{\mbox{$C^{1}_{0}(#1)$}}
\newcommand{\lspac}[2]{\mbox{$L^{#1}(#2)$}}
\newc{\hspac}[2]{\mbox{$C^{0,#1}(#2)$}}
\newc{\Hspac}[2]{\mbox{$C^{1,#1}(#2)$}}
\newc{\Hosp}{\mbox{$H^{1}_{0}$}}
\newcommand{\Lsp}[1]{\mbox{$L^{#1}(\Og)$}}
\newc{\hsp}{\Hosp(\Og)}


\newc{\ag}{\alpha}
\newc{\bg}{\beta}
\newc{\cg}{\gamma}\newc{\Cg}{\Gamma}
\newc{\dg}{\delta}\newc{\Dg}{\Delta}
\newc{\eg}{\varepsilon}
\newc{\zg}{\zeta}
\newc{\thg}{\theta}
\newc{\llg}{\lambda}\newc{\LLg}{\Lambda}
\newc{\kg}{\kappa}
\newc{\rg}{\rho}
\newc{\sg}{\sigma}\newc{\Sg}{\Sigma}
\newc{\tg}{\tau}
\newc{\fg}{\phi}\newc{\Fg}{\Phi}
\newc{\vfg}{\varphi}
\newc{\og}{\omega}\newc{\Og}{\Omega}
\newc{\pdh}{\partial}

\newc{\ccG}{{\cal G}}


\newc{\ii}[1]{\int_{#1}}
\newc{\iidx}[2]{{\dspl\int_{#1}~#2~dx}}
\newc{\bii}[1]{{\dspl \ii{#1} }}
\newc{\biii}[2]{{\dspl \iii{#1}{#2} }}
\newc{\su}[2]{\sum_{#1}^{#2}}
\newc{\bsu}[2]{{\dspl \su{#1}{#2} }}

\newcommand{\iiomdx}[1]{{\dspl\int_{\Og}~ #1 ~dx}}
\newc{\biiom}[1]{{\dspl\int_{\bdrom}~ #1 ~d\sg}}
\newc{\io}[1]{{\dspl\int_{\Og}~ #1 ~dx}}
\newc{\bio}[1]{{\dspl\int_{\bdrom}~ #1 ~d\sg}}
\newc{\bsir}{\bsu{i=1}{r}}
\newc{\bsim}{\bsu{i=1}{m}}

\newc{\iibr}[2]{\iidx{\bprw{#1}}{#2}}
\newc{\Intbr}[1]{\iibr{R}{#1}}
\newc{\intbr}[1]{\iibr{\rg}{#1}}
\newc{\intt}[3]{\int_{#1}^{#2}\int_\Og~#3~dxdt}

\newc{\itQ}[2]{\dspl{\int\hspace{-2.5mm}\int_{#1}~#2~dz}}
\newc{\mitQ}[2]{\dspl{\rule[1mm]{4mm}{.3mm}\hspace{-5.3mm}\int\hspace{-2.5mm}\int_{#1}~#2~dz}}
\newc{\mitQQ}[3]{\dspl{\rule[1mm]{4mm}{.3mm}\hspace{-5.3mm}\int\hspace{-2.5mm}\int_{#1}~#2~#3}}

\newc{\mitx}[2]{\dspl{\rule[1mm]{3mm}{.3mm}\hspace{-4mm}\int_{#1}~#2~dx}}

\newc{\mitQq}[2]{\dspl{\rule[1mm]{4mm}{.3mm}\hspace{-5.3mm}\int\hspace{-2.5mm}\int_{#1}~#2~d\bar{z}}}
\newc{\itQq}[2]{\dspl{\int\hspace{-2.5mm}\int_{#1}~#2~d\bar{z}}}

\newc{\pder}[2]{\dspl{\frac{\partial #1}{\partial #2}}}


\newc{\ui}{u_{i}}
\newcommand{\upl}{u^{+}}
\newcommand{\umn}{u^{-}}
\newcommand{\un}{\{ u_{n}\}}

\newcommand{\uo}{u_{0}}
\newc{\voi}{v_{i}^{0}}
\newc{\uoi}{u_{i}^{0}}
\newc{\vu}{\vec{u}}

\newc{\xo}{x_{0}}
\newc{\Br}{B_{R}}
\newc{\Bro}{\Br (\xo)}
\newc{\bdrom}{\bdry{\Og}}
\newc{\ogr}[1]{\Og_{#1}}
\newc{\Bxo}{B_{x_0}}

\newc{\inP}[2]{\|#1(\bullet,t)\|_#2\in\cP}
\newc{\cO}{{\mathcal O}}
\newc{\inO}[2]{\|#1(\bullet,t)\|_#2\in\cO}

\newc{\newl}{\\ &&}

\newc{\bilhom}{\mbox{Bil}(\mbox{Hom}(\RR^{nm},\RR^{nm}))}
\newc{\VV}[1]{{V(Q_{#1})}}

\newc{\ccA}{{\mathcal A}}
\newc{\ccB}{{\mathcal B}}
\newc{\ccC}{{\mathcal C}}
\newc{\ccD}{{\mathcal D}}
\newc{\ccE}{{\mathcal E}}
\newc{\ccH}{\mathcal{H}}
\newc{\ccF}{\mathcal{F}}
\newc{\ccI}{{\mathcal I}}
\newc{\ccJ}{{\mathcal J}}
\newc{\ccP}{{\mathcal P}}
\newc{\ccQ}{{\mathcal Q}}
\newc{\ccR}{{\mathcal R}}
\newc{\ccS}{{\mathcal S}}
\newc{\ccT}{{\mathcal T}}
\newc{\ccX}{{\mathcal X}}
\newc{\ccY}{{\mathcal Y}}
\newc{\ccZ}{{\mathcal Z}}

\newc{\bb}[1]{{\mathbf #1}}
\newc{\bbA}{{\mathbf A}}
\newc{\myprod}[1]{\langle #1 \rangle}
\newc{\mypar}[1]{\left( #1 \right)}

\newc{\BLLg}{\mathbf{\LLg}}

\newc{\mA}{\mathbf{A}}
\newc{\mB}{\mathbf{B}}
\newc{\mC}{\mathbf{C}}
\newc{\mD}{\mathbf{D}}
\newc{\mE}{\mathbf{E}}
\newc{\mF}{\mathbf{F}}
\newc{\mJ}{\mathbf{J}}
\newc{\mG}{\mathbf{G}}
\newc{\mP}{\mathbf{P}}
\newc{\mR}{\mathbf{R}}
\newc{\mQ}{\mathbf{Q}}
\newc{\mX}{\mathbf{X}}
\newc{\muu}{\mathbf{u}}
\newc{\mvv}{\mathbf{v}}

\newc{\mllg}{\mathbb{\lambda}}
\newc{\mLLg}{\mathbf{\LLg}}


\newc{\lspn}[2]{\mbox{$\| #1\|_{\Lsp{#2}}$}}
\newc{\Lpn}[2]{\mbox{$\| #1\|_{#2}$}}
\newc{\Hn}[1]{\mbox{$\| #1\|_{H^1(\Og)}$}}


\newc{\cyl}[1]{\og\times \{#1\}}
\newc{\cyll}{\og\times[0,1]}
\newc{\vx}[1]{v\cdot #1}
\newc{\vtx}[1]{v(t,x)\cdot #1}
\newc{\vn}{\vx{n}}

\newcommand{\RR}{{\rm I\kern -1.6pt{\rm R}}}


\newenvironment{proof}{\noindent\textbf{Proof.}\ }
{\nopagebreak\hbox{ }\hfill$\Box$\bigskip}


\newc{\itQQ}[2]{\dspl{\int_{#1}#2\,dz}}
\newc{\mmitQQ}[2]{\dspl{\rule[1mm]{4mm}{.3mm}\hspace{-4.3mm}\int_{#1}~#2~dz}}
\newc{\MmitQQ}[2]{\dspl{\rule[1mm]{4mm}{.3mm}\hspace{-4.3mm}\int_{#1}~#2~d\mu}}

\newc{\MUmitQQ}[3]{\dspl{\rule[1mm]{4mm}{.3mm}\hspace{-4.3mm}\int_{#1}~#2~d#3}}
\newc{\MUitQQ}[3]{\dspl{\int_{#1}~#2~d#3}}

\vspace*{-.8in}
\begin{center} {\LARGE\em Local and Global Existence of Strong Solutions to Large Cross Diffusion Systems.}

 \end{center}

\vspace{.1in}

\begin{center}

{\sc Dung Le}{\footnote {Department of Mathematics, University of
Texas at San
Antonio, One UTSA Circle, San Antonio, TX 78249. {\tt Email: Dung.Le@utsa.edu}
}}

\end{center}

\begin{abstract}
We study the solvability of a general class of cross diffusion systems and establish the local and global existence of their strong solutions under the weakest assumption that they are VMO. This work simplifies the setting in our previous work \cite{dleANS} and provides new extensions which are more verifiable in applications. \end{abstract}

\vspace{.2in}

\section{Introduction}\label{introsec}\eqnoset

In this paper, for any $T_0>0$ and bounded domain $\Og$ with smooth boundary in $\RR^n$, $n\ge2$, we consider the following general parabolic system of $m$ equations ($m\ge2$) 
\beqno{e1}\left\{\barr{ll} u_t=\Div(A(u)Du)+\hat{f}(u,Du)& (x,t)\in Q=\Og\times(0,T_0),\\u(x,0)=U_0(x)& x\in\Og\\\mbox{$u=0$ or $\frac{\partial u}{\partial \nu}=0$ on $\partial \Og\times(0,T_0)$}. &\earr\right.\eeq
where $A(u)$ is a $m\times m$ matrix in $u$; $u:\Og\to\RR^m$, $\hat{f}:\RR^m\times \RR^{mn}\times\RR^{nm}\to\RR^m$ are vector valued functions. 
The initial data $U_0$ is given in $W^{1,p_0}(\Og,\RR^m)$ for some $p_0>n$, the dimension of $\Og$. As usual, $W^{1,p}(\Og,\RR^m)$, $p\ge1$, will denote the standard Sobolev spaces whose elements are vector valued functions $u\,:\,\Og\to \RR^m$ with finite norm $$\|u\|_{W^{1,p}(\Og,\RR^m)} = \|u\|_{L^p(\Og)} + \|Du\|_{L^p(\Og)}.$$ 

We say that $u$ is a strong solution if $u$ is continuous on $\bar{Q}$ with $Du\in L^\infty_{loc}(Q)$ and $D^2u\in L^2_{loc}(Q)$.

The strongly coupled system \mref{e1} appears in many physical applications, for instance, Maxwell-Stephan systems describing the diffusive transport of multicomponent mixtures, models in reaction and diffusion in electrolysis, flows in porous media, diffusion of polymers, or population dynamics, among others. We refer the reader to the recent work \cite{JA} and the references therein for the models and the existence of their {\em weak} solutions. 

The first fundamental problem in the study of \mref{e1} is the local and global existence of its solutions. One can decide to work with either weak or strong solutions. In the first case, the existence of a weak solution can be easily achieved via Galerkin or variational methods but its regularity (e.g., boundedness, H\"older continuity of the solution and its higher derivatives) is still an open issue. Several works have been done along this line to improve the early work \cite{GiaS} and establish {\em partial regularity} of {\em bounded} weak solutions to \mref{e1}.

Otherwise, if strong solutions are considered then theirs existence can be established via semigroup theories as in the works of Amann \cite{Am1,Am2}. Combining with interpolation theories of Sobolev's spaces, Amann established local and global existence of a strong solution $u$ of \mref{e1} under the assumption that one can controll $\|u\|_{W^{1,p}(\Og,\RR^m)}$ for some $p>n$.

In both forementioned approaches, the assumption on the boundedness of $u$ must be the starting point. 
For strongly coupled systems like \mref{e1}, as invariant/maximum principles for cross diffusion systems are generally unavailable, the boundedness of the solutions is already a hard problem. One usually needs to use ad hoc techniques on the case by case basis to show that $u$ is bounded (see \cite{kuf1,Red}).  Even for bounded weak solutions, we know that they are only H\"older continuous almost everywhere (see \cite{GiaS}).
In addition, there are counter examples for systems ($m>1$) which exhibit solutions that start smoothly and remain bounded but develop singularities in higher norms in finite times (see \cite{JS}).

In our recent work \cite{dleANS}, we choose a different approach making use of fixed point theory and discuss the solvability of \mref{e1} under the weakest assumption that $u$ is VMO (see \mref{KeyVMO} below) and much more general structural conditions on the data of \mref{e1}.  The proof in \cite{dleANS} relies on fixed point theories, instead of the semigroup approach in \cite{Am1}, and weighted Gagliardo-Ninrehnberg inequalities involving BMO norms.

In particular, we assumed in \cite{dleANS} the followings.

\bdes

\item[A)] $A(u)$ is $C^1$ in $u$ and there are constants $\llg_0,C_*>0$  and a scalar $C^1$ function $\llg(u)$  such that  for all $u\in\RR^m$ and $\zeta\in\RR^{nm}$ 
\beqno{A1} \llg(u)\ge \llg_0,\;\llg(u)|\zeta|^2 \le \myprod{A(u)\zeta,\zeta} \mbox{ and } |A(u)|\le C_*\llg(u).\eeq 

In addition, $|A_u|\le C|\llg_u|$ and  the following number  is finite:\beqno{LLgx} \mathbf{\LLg}=\sup_{W\in\RR^m}\frac{|\llg_W(W)|}{\llg(W)}.\eeq

\edes 

Here and throughout this paper, if $B$ is a $C^1$ function in $u\in \RR^m$ then we abbreviate it  derivative $\frac{\partial B}{\partial u}$ by $B_u$.

\bdes \item[F)] There exist a constant $C$ and a $C^1$ function $f:\RR^m\to\RR^m$ such that for any diffrentiable vector valued functions $u:\RR^n\to\RR^m$ and $p:\RR^n\to\RR^{mn}$ \beqno{FUDU11}|\hat{f}(u,p)| \le C\llg^\frac12(u)|p| + f(u),\eeq 
\beqno{FUDU11a}|D\hat{f}(u,p)| \le C\llg^\frac12(u)|Dp| + C\frac{|\llg_u(u)|}{\llg^\frac12(u)}|Du||p|+|f_u(u)||Du|,\eeq
\beqno{fUDU21}|f_u(u)| \le C\llg(u).\eeq
\edes

The local existence of a strong solution of \mref{e1} was proved in \cite{dleANS} under the key assumption (see the condition M') in the paper) that any strong solution $u$ of the system satisfies:  
for any given $\mu_0>0$ there is a positive $R_{\mu_0}$ such that   \beqno{KeyVMO} \mathbf{\LLg}^2\sup_{x_0\in\bar{\Og},t\in (0,T_0)}\|u(\cdot,t)\|_{BMO(B_{R_{\mu_0}}(x_0)\cap\Og)}^2 \le \mu_0.\eeq  

Here, a locally integrable function $U:\Og\to\RR^m$ is said to be in $BMO(\Og)$ if the quantity $$[U]_*:=\sup_{B_R(y)\subset\Og} \mitx{B_R(y)}{|U-U_{B_R(y)}|}\quad\mbox{is finite}.$$

The Banach space $BMO(\Og,\RR^m)$ consists of functions with finite norm $$\|U\|_{BMO(\Og,\RR^m)}:=[U]_*+\|U\|_{L^1(\Og,\RR^m)}.$$
We also say that $U$ is VMO on $\Og$ if $\inf_{R>0,B_R\subset\Og}\|U\|_{BMO(B_R,\RR^m)}=0$.

In this paper, for simplicity of the presentation and with models in applications in mind, we consider first the following special form of the reaction terms which are linear in $Du$, namely,
$\hat{f}(u,Du) = B(u)Du + f(u)$ and study local and global existence of strong solutions. Thanks to this form of $\hat{f}$ the fixed point argument in \cite{dleANS} can be greatly simplified. Furthermore, we will provide conditions which are a bit stronger than M') but verifiable in applications. 
In particular, we will show that a strong solution $u$ exists globally if the norm $\|u\|_{W^{1,n}(\Og)}$ does not blow up in finite time. This relaxes Amann's conditions in \cite{Am1} which required a control on $\|u\|_{W^{1,p}(\Og)}$ for some $p>n$.  Again, we are not assuming that $u$ is bounded and our structural conditions A) and F) are more general than those in \cite{Am1}.
We then show that our results continue to hold for general $\hat{f}(u,Du)$ with linear growth in $Du$.

We organize our paper as follows. In \refsec{res} we state our main results. In \refsec{GNsec} we give another version of the local weighted Gagliardo-Nirenberg inequality \cite[Lemma 2.4]{dleANS}, which is one of the main ingredients of the proof in \cite{dleANS} and our main theorem in this paper. Technical results and auxiliary lemmas needed for tthe proof of the main results for linear reaction terms will be given in \refsec{w12est}. The proof of the general case is given in \refsec{genreaction}.

\section{Preliminaries and Main Results}\eqnoset\label{res}

We state the main results of this paper in this this section. Our first main result concerns the local existence of strong solutions to \mref{e1} with $\hat{f}$ being linear in $Du$.
\beqno{fspecial} \hat{f}(u,Du) = B(u)Du+f(u).\eeq 

We imbed \mref{e1} in the following family of systems 
\beqno{e1famzzz}\left\{\barr{l} u_t=\Div(A(\sg u)Du)+\hat{f}(\sg u,\sg Du)\quad (x,t)\in Q=\Og\times(0,T_0), \sg\in[0,1]\\u(x,0)=U_0(x)\quad  x\in\Og\\\mbox{$u=0$ or $\frac{\partial u}{\partial \nu}=0$ on $\partial \Og\times(0,T_0)$}. \earr\right.\eeq

In \cite{dleANS} we also assumed the condition SG) which requires that the eigenvalues of the matrix $A(u)$ are not too far apart. Namely, we need that $(n-2)/n <C_*^{-1}$, where $C_*$ is, in certain sense, the ratio of the largest and smallest eiegenvalues of $A(u)$.
Being inspired by this condition, 
We then assume that there is $n^*>n$ be such that $(n^*-2)/n^* = C_*^{-1}$ or
\beqno{nSGdef} n^*= \frac{C_*}{2(C_*-1)}.\eeq

\btheo{dleNL-mainthm} 
Assume A), F) with \mref{fspecial}. There is $\mu_0>0$ sufficiently small in terms of the constants in A) and F) such that if the following two conditions hold for {\em any strong} solution $u$ of \mref{e1famzzz} 

\bdes \item[M.1)]
For some positive $R_{\mu_0}$, which may depends on $T_0$,  \beqno{Keymu0} \mathbf{\LLg}^2\sup_{x_0\in\bar{\Og},t\in (0,T_0)}\|u(\cdot,t)\|_{BMO(B_{R}(x_0)\cap\Og)}^2 \le \mu_0;\eeq  

\item[M.2)] The following quantity \beqno{keydupANSp2zz}
C_{T_0}:=\itQ{\Og\times(0,T_0)}{|Du|^{2}} \mbox{ is finite};\eeq

\item[L)] There are constants $L(T_0)$ and $r^*>\frac{n}{p^*-n}$, with $p^*=\min\{n^*,p_0\}$, see the definition of $n^*$ in \mref{nSGdef}, such that
\beqno{llghyp}\sup_{t\in(0,T_0)}\|\llg(u(\cdot,t))\|_{L^{r^*}(\Og)}\le L(T_0),\eeq \edes
then \mref{e1} has a unique strong solution $u$ on $\Og\times(0,T_0)$.
\etheo

As an immediate consequence of this theorem and \reflemm{Uniquelem} at the end of the paper, we have the following result on the global existence of strong solutions.

\bcoro{dleNL-maincoro} If the assumptions of \reftheo{dleNL-mainthm} hold for all $T_0>0$ then \mref{e1} has a unique strong solution $u$ which exists globally on $\Og\times(0,\infty)$. \ecoro

The next results are more applicable and improve those of Amann in \cite{Am1,Am2}. Basically, we need only to controll the $W^{1,n}(\Og)$ norm of strong solutions while \cite{Am1,Am2} required that their $W^{1,p}(\Og)$ norms do not blow up in finite time for some $p>n$, and thus the boundedness of the solutions is needed.

\bcoro{dlec1coro} The conclusion of \reftheo{dleNL-mainthm} holds if M.1) and M.2) are replaced by the following assumption.
\bdes \item[D)] There is a constant $C_{T_0}$ such that for any $t\in(0,T_0)$
\beqno{pis20az5} \|u(\cdot,t)\|_{W^{1,n}(\Og)}\le C_{T_0}.\eeq \edes

If this condition holds for all $T_0>0$ then $u$ exists globally.
\ecoro

Finally, concerning the integrability condition of $\llg(u)$ in L), we have the following result.

\bcoro{dlec2coro} If
we assume further that there are constants $\BLLg_1,\eg_0>0$ such that \beqno{llgullgeg} |\llg_u(u)|\le \BLLg_1\llg^{1-\eg_0}(u)\quad \forall u\in\RR^m.\eeq 

Then the conclusion of \refcoro{dlec1coro} holds if L) is replaced by the following weaker one.

\bdes\item[L')] There are constants $L(T_0)$  $s_0>0$ \beqno{llgs0zzz5}\sup_{t\in(0,T_0)}\|\llg^{s_0}(u(\cdot,t))\|_{L^{1}(\Og)}\le L(T_0). \eeq \edes
\ecoro

It is easy to see that the condition \mref{llgullgeg} holds if $\llg(u)$ with polynomial growth in $u$.

We now turn to the general \mref{e1} and consider general reaction terms. We imbed \mref{e1} in the following family of systems 
\beqno{e1famza}\left\{\barr{l} u_t=\Div(A(\sg u)Du)+\hat{f}(\sg u,\sg Du)\quad (x,t)\in Q=\Og\times(0,T_0), \sg\in[0,1]\\u(x,0)=U_0(x)\quad  x\in\Og\\\mbox{$u=0$ or $\frac{\partial u}{\partial \nu}=0$ on $\partial \Og\times(0,T_0)$}. \earr\right.\eeq

\btheo{Genthm} Assume A) and F). We assume further that
\beqno{DFDu0} |\partial_\zeta \hat{f}(u,\zeta)|\le C\llg^\frac12(u),\; |\partial_u \hat{f}(u,0)|\le C\llg(u)\quad \forall u\in\RR^m,\forall\zeta\in\RR^{nm}.\eeq

Then the conclusions of \reftheo{dleNL-mainthm}, \refcoro{dlec1coro} and \refcoro{dlec2coro} continue to hold under their asumptions M.1)-L), D) and L') respectively. \etheo

\section{Technical results} \eqnoset\label{GNsec}
We provide the following version of the local weighted Gagliardo-Nirenberg inequality \cite[Lemma 2.4]{dleANS}, which is one of the main ingredients of the proof in \cite{dleANS} and our main theorem in this paper. In order to state the assumption for this type of inequalities, we recall some well known notions from Harmonic Analysis. For $\cg\in(1,\infty)$ we say that a nonnegative locally integrable function $w$ is an $A_\cg$ weight if the quantity
\beqno{aweight} [w]_{\cg} := \sup_{B_R(y)\subset\Og} \left(\mitx{B_R(y)}{w}\right) \left(\mitx{B_R(y)}{w^{1-\cg'}}\right)^{\cg-1} \quad\mbox{is finite}.\eeq
Here, $\cg'=\cg/(\cg-1)$. For more details on these classes we refer the reader to \cite{OP,st}. 

We also make use of Hardy spaces ${\cal H}^{1}$. For any $y\in\Og$ and $\eg>0$, let $\fg$ be any function in $C^\infty_0(B_1(y))$ with $|D\fg|\le C_1$. Let $\fg_\eg(x)=\eg^{-n}\fg(\frac{x}{\eg})$ (then $|D\fg_\eg|\le C_1\eg^{-1-n}$). From \cite{st}, a function $g$ is in ${\cal H}^{1}(\Og)$ if \beqno{Hardynorm}\sup_{\eg>0}g*\fg_\eg\in L^1(\Og) \mbox{ and } \|g\|_{{\cal H}^1} = \|g\|_{L^1(\Og)}+\|\sup_{\eg>0}g*\fg_\eg\|_{L^1(\Og)}.\eeq

Throughout this paper,  when there is no ambiguity $C, C_i$ will denote universal constants that can change from line to line in our argument. If necessary, $C(\cdots)$ or $C_{(\cdots)}$ are used to denote quantities which are bounded in terms of theirs parameters in $(\cdots)$. We will also write $a\sim b$ if there are two generic positive constants $C_1,C_2$ such that $C_1b \le a \le C_2b$. Furthermore, we denote by $B_R(x_0)$ a ball with center $x_0\in\bar{\Og}$. In the sequel, if the center $x_0$ is already specified then we simply write $B_R,\Og_R$ for $B_R(x_0),B_R(x_0)\cap\Og$ respectively. 

We have the following version of \cite[Lemma 2.4]{dleANS}.

\blemm{GNlocalnew} Let $u,U:\Og\to \RR^m$ be vector-valued functions with $u\in C^1(\Og)$, $U\in C^2(\Og)$, and let $\Fg:\RR^m\to\RR$ be a $C^1$ function and $\Fg(u)^\frac{2}{p+2}$ be an $A_{p/(p+2)+1}$ weight. Suppose that either $U$ or $\Fg^2(u)\frac{\partial U}{\partial \nu}$ vanish on the boundary $\partial \Og$  of $\Og$. For any ball $B_t(x_0)$ with center $x_0\in\bar{\Og}$
we set \beqno{Ideft} I_1(t):=\iidx{\Og_t}{\Fg^2(u)|DU|^{2p+2}},\,\hat{I}_1(t):=\iidx{\Og_t}{\Fg^2(u)|Du|^{2p+2}},\eeq
\beqno{Idef1t} \bar{I}_1(t):=\iidx{\Og_t}{|\Fg_u(u)|^2(|DU|^{2p+2}+|Du|^{2p+2})},\eeq and \beqno{Idef2t}
I_2(t):=\iidx{\Og_t}{\Fg^2(u)|DU|^{2p-2}|D^2U|^2}.\eeq

Consider any ball $B_s$ concentric with $B_t$, $0<s<t$, and any nonnegative $C^1$ function $\psi$  such that $\psi=1$ in $B_s$ and $\psi=0$ outside $B_t$. Then, for any $\eg>0$ there are positive constants $C_{\eg,\Fg}$, which depends on $[\Fg]_{p/(p+2)+1}$, and  $C_\eg$ such that \beqno{GNlocalest}\barr{lll}I_1(s)&\le& \eg[I_1(t)+\hat{I}_1(t)]+C_{\eg,\Fg}\|U\|^2_{BMO(\Og_t)}
\left[\bar{I}_1(t) +  I_2(t)\right]\\ &&+C_\eg\|U\|_{BMO(\Og_t)}^2\sup_{x\in B_t}|D\psi(x)|^2\iidx{\Og_t}{\Fg^{2}(u)|DU|^{2p}}.\earr\eeq
\elemm

The only difference between the two versions is that the factor $\|U\|_{BMO(B_t)}^2$  in the last terms of  \mref{GNlocalest} here replaces the factor $\|U\|_{BMO(B_t)}$  in the eqn. (2.17) of \cite[Lemma 2.4]{dleANS}. The two proofs differ only by the order of using Young's inequality in the argument. For the sake of completeness we present the details.

\bproof We revisit the proof of the global weighted Gagliardo-Nirenberg inequality \cite[Lemma 2.1]{dleANS}. Integrating by parts and using the assumptions on $\Fg(u),U$ on $\partial\Og$, we have
\beqno{FSstart}\iidx{\Og}{\Fg^2(u)\psi^2|DU|^{2p+2}} = -\iidx{\Og}{U\Div(\Fg^2(u)\psi^2|DU|^{2p}DU)}.\eeq

We will show that  $g=\Div(\Fg^2\psi^2|DU|^{2p}DU)$ belongs to the Hardy space ${\cal H}^1$ (see \mref{Hardynorm}).
We write $g=g_1+g_2$ with $g_i=\Div V_i$, setting  
$$V_1= \Fg(u) \psi|DU|^{p+1}\left(\Fg(u)\psi |DU|^{p-1}DU-\mitx{\Og_\eg}{\Fg(u)\psi |DU|^{p-1}DU}\right),$$ (here, $\Og_\eg=B_\eg(y)\cap\Og$, with $y$ being the variable of $V_i$) and
$$V_2= \Fg(u)\psi |DU|^{p+1}\mitx{\Og_\eg}{\Fg(u)\psi |DU|^{p-1}DU}.$$

In estimating $V_1$ we follow the proof of \cite[Lemma 2.1]{dleANS} and replace $\Fg(u)$  by $\Fg(u)\psi(x)$. There will be some extra terms in the proof in computing $D(\Fg(u)\psi)$. In particular, in estimating  $Dh$ in the right hand side of \cite[eqn. (2.8)]{dleANS} we have the following term and it can be estimated as follows $$\left(\mitx{\Og_\eg}{\Fg^{s_*}(u)|D\psi|^{s_*}|DU|^{ps_*}}\right)^\frac1{s_*} \le \sup_{x\in B_t}|D\psi|\left(\mitx{\Og_\eg}{\Fg^{s_*}(u)|DU|^{ps_*}}\right)^\frac1{s_*}.$$

We then use  H\"older's inequality, in the right hand side of \cite[eqn. (2.9)]{dleANS} (with $\Og=B_t$, $M$ is the maximal operator, and the definition of $\Psi_1$ there)
$$\barr{ll}\lefteqn{\sup|D\psi|\iidx{\Og_\eg}{\Psi_1M(\Fg^{s_*}(u)|DU|^{ps_*})^\frac{1}{s_*}} \le}\hspace{2cm}&\\& C \sup_{x\in B_t}|D\psi|\left(\iidx{\Og_\eg}{\Psi_1^2}\right)^\frac12\left(\iidx{\Og_t}{M(\Fg^{s_*}(u)|DU|^{ps_*})^\frac{2}{s_*}}\right)^\frac12.\earr$$

The last integral can be bounded via the Hardy-Littlewood inequality \cite[eqn. (2.4)]{dleANS} by $$\iidx{\Og_\eg}{\Fg^{2}(u)|DU|^{2p}}.$$

Using the fact that $|\psi|\le 1$ and $\Og=B_t$, the proof of \cite[Lemma 2.1]{dleANS} can go on and \cite[eqn. (2.11)]{dleANS} now becomes \beqno{g1estz}\barr{ll}\lefteqn{\iidx{\Og_t}{\sup_\eg|g_1*\fg_\eg|} \le}\hspace{1cm} &\\&C
\left[I_1^\frac12\bar{I}_1^\frac12+I_1^\frac12 I_2^\frac12 +\sup_{B_t}|D\psi|I_1^\frac12\left(\iidx{\Og_t}{\Fg^{2}(u)|DU|^{2p}}\right)^\frac12\right].\earr\eeq

Similarly, in considering $g_2=\mbox{div}V_2$, we will have the following extra term $$\sup_\eg |\fg_\eg*\Fg(u)||D\psi||DU|^{p+1}\mitx{\Og_\eg}{\Fg(u)|DU|^p},$$ which can be estimated by $ \sup_{B_t}|D\psi|M(\Fg(u)|DU|^{p+1})M(\Fg(u)|DU|^p)$. Using H\"older's inequality and the Hardy-Littlewood inequality, the integral over $B_t$ of this quantity is bounded by
$$ \sup_{B_t}|D\psi|I_1^\frac12\left(\iidx{\Og_\eg}{\Fg^{2}(u)|DU|^{2p}}\right)^\frac12.$$

Therefore the estimate \cite[eqn. (2.16)]{dleANS} is now \mref{g1estz} with $g_1$ being replaced by $g_2$. Combining the estimates for $g_1,g_2$ and using Young's inequality, we get  $$\barr{ll}\lefteqn{\iidx{\Og_t}{\sup_\eg|g*\fg_\eg|}\le}\hspace{.5cm} &\\&\sup_{B_t}|D\psi|I_1^\frac12\left(\iidx{\Og_t}{\Fg^{2}(u)|DU|^{2p}}\right)^\frac12 +C_\Fg
\left[\bar{I}_1^\frac12(I_1^\frac12+\hat{I}_1^\frac12) + I_1^\frac12 I_2^\frac12\right].\earr$$

The above gives an estimate for the ${\cal H}^1$ norm of $g$ (see \mref{Hardynorm}). By the Fefferman-Stein duality theorem (see \cite{st}) and \mref{FSstart}, we obtain $$\barr{ll}\lefteqn{\iidx{\Og_t}{\Fg^2(u)\psi^2|DU|^{2p+2}}\le \|U\|_{BMO(\Og_t)}\|g\|_{{\cal H}^1} \le}\hspace{.5cm} &\\&\|U\|_{BMO(\Og_t)}\left(\sup_{B_t}|D\psi|I_1^\frac12\left(\iidx{\Og_t}{\Fg^{2}(u)|DU|^{2p}}\right)^\frac12 +C_\Fg
\left[\bar{I}_1^\frac12(I_1^\frac12+\hat{I}_1^\frac12) + I_1^\frac12 I_2^\frac12\right]\right) .\earr$$

A simple use Young's inequality and then the fact that $\psi=1$ in $B_s$ gives \mref{GNlocalest} and completes the proof. \eproof
	
	We now let $\Fg\equiv1$ and $\psi$ be a cutoff function for $B_s,B_t$, i.e. $\psi=1$ in $B_s$ and $\psi=0$ outside $B_t$ and $|D\psi|\le 1/(t-s)$, then $\Fg$ is an $A_\cg$ weight for all $\cg>1$ and $\Fg_u\equiv0$. The following version of the above lemma with $u=U$ suffices for our purpose in this paper.
	\blemm{GNlocalspec} Let $U:\Og\to \RR^m$ be a vector-valued function in $C^2(\Og)$. Suppose that either $U$ or $\frac{\partial U}{\partial \nu}$ vanish on the boundary $\partial \Og$  of $\Og$.
	For any two concentric balls $B_s,B_t$, with $s<t$, and any $p\ge1,\eg>0$ there is $C_\eg>0$ such that
	\beqno{GNlocalestspec}\barr{ll}\lefteqn{\iidx{\Og_s}{|DU|^{2p+2}}\le \eg\iidx{B_t}{|DU|^{2p+2}}+}\hspace{2cm}&\\& C_{\eg}\|U\|^2_{BMO(\Og_t)}\iidx{\Og_t}{[|DU|^{2p-2}|D^2U|^2+(t-s)^{-2}|DU|^{2p}]}.\earr\eeq
	
	\elemm

\section{The proof of the main results}\eqnoset\label{w12est}

In this section, we present the proof of \reftheo{dleNL-mainthm} and its corollaries. The proof relies on the Leray Schauder theorem. We obtain the existence of a strong solution $u$ of \mref{e1} as a fixed point of a nonlinear map defined on an appropriate Banach space.

Let us consider the Banach space $\mX=C(Q,\RR^m)$, where $Q=\Og\times(0,T_0)$.
For any given $u\in\mX$ and $\sg\in[0,1]$, we consider the following linear system
\beqno{Tmapdef}\left\{\barr{l} w_t=\Div(A(\sg u)Dw)+B(\sg u)Dw+f(\sg u)\quad (x,t)\in Q, \\w(x,0)=U_0(x)\quad  x\in\Og\\\mbox{$w=0$ or $\frac{\partial w}{\partial \nu}=0$ on $\partial \Og\times(0,T_0)$}. \earr\right.\eeq

We then define $T_\sg(u)=w$. It is clear that a fixed point of $T_\sg$ solves \mref{e1famzzz}. In order to apply the Leray-Schauder theorem, we need to establish the followings.

\bdes \item[Step 1]The map $T_\sg:\mX\to\mX$ is well defined and compact. 
\item[Step 2] There is a constant $M$ such that $\|u\|_\mX\le M$ for any fixed points of $u=T_\sg(u)$.
\edes

The checking of Step 1 is fairly standard thanks to the following lemma.

\blemm{Tmaplem} The map $T_\sg:\mX\to\mX$ is well defined and compact. \elemm

\bproof For each $u\in\mX$, $A(\sg u)$ satisfies the ellipticity condition A) and the data of the linear system \mref{Tmapdef} are bounded and continuous. So that \mref{Tmapdef} satisfies the assumptions of Theorem 1.1 in \cite[Chapter VII]{LSU}, which applies to the system
\beqno{LSUsys} w_t=\Div(\mathbf{a} Dw) + \mathbf{b}Dw+\mathbf{g},\eeq under the assumption that $\mathbf{a},\mathbf{b}$ and $\|\mathbf{g}\|_{q,r,Q}$ are bounded for sufficiently large $q,r$ such that $1/r+n/(2q)=1$. Here, we denoted for any vector valued function $F$
$$\|F\|_{q,r,Q}:=\left(\int_0^{T_0}\left(\iidx{\Og}{|F(x,t)|^q}\right)^\frac{r}{q}dt\right)^\frac{1}{r}.$$

Theorem 1.1 in \cite[Chapter VII]{LSU} shows that $w$ exists uniquely so that $T_\sg(u)$ is well defined. Moreover, as the initial condition $w(\cdot,0)=U_0(x)$ belongs to $W^{1,p_0}(\Og)$ and then $C^\bg_0(\Og)$ for $\bg_0=1-n/p_0>0$, a combination of Theorems 2.1 and 3.1 in \cite[Chapter VII]{LSU} shows that $w$ belongs to $C^{\ag_0,\ag_0/2}(\bar{Q},\RR^m)$ for some $\ag_0>0$. Moreover, the norm $\|w\|_{C^{\ag_0,\ag_0/2}(\bar{Q})}$ depending on $\bg_0$ and $\|A(\sg u)\|_\infty$, $\|B(\sg u)\|_\infty$ $\|f_\sg(u)\|_{q,r,Q}$. 
Thus, if $u$ belongs to a bounded set $K$ of $\mX$ then $\|u\|_\mX\le M$ for some $M$ and there is a constant $C$ such that $$\|w\|_{C^{\ag_0,\ag_0/2}(\bar{Q})}\le C(M,\|U_0(\cdot,0)\|_{C^{\bg_0}(\Og)}).$$
Hence, $T_\sg(K)$ is  compact in $\mX$ and
$T_\sg:\mX\to\mX$ is a compact map. \eproof

We now turn to Step 2, the hardest part of the proof, and provide a uniform estimate for the fixed points of $T_\sg$. Such a fixed point $u$ of $T_\sg$  satisfies \mref{Tmapdef} and belongs to $\mX$. Therefore, $u$ is a bounded weak solution and continuous so that \cite[Theorems 2.1 and 3.2]{GiaS}  apply and yield that $Du$ is bounded in $\Og\times(t_0,T_0)$ for all $t_0>0$. Thus, $Du$ is locally bounded in $\Og\times(0,T_0)$. It is then well known that $D^2u$ exists in $L^2_{loc}(\Og\times(0,T_0)$ and $u$ is a  strong solution in  $\Og\times(0,T_0)$. 

Thus, in the rest of this section, we consider a strong solution $u$ of \mref{Tmapdef}. As the data of \mref{Tmapdef} satisfy the structural conditions A), F) with the same set of constants and the asumptions M.1)-L) are assumed to be uniform for all $\sg\in[0,1]$, we will only present the proof for $\sg=1$ in the sequel.

We should also emphasize that the estimates in the rest of this section do not require the special form of $\hat{f}$ in \mref{fspecial} but the growth condition in F).

For any two concentric balls $B_s,B_t$ with $s<t$, we say that $\psi$ is a cutoff function for $B_s,B_{t}$ if $\psi$ is a $C^1$ function satisfying $\psi\equiv1$ in $B_s$ and $\psi\equiv0$ outside $B_t$ and $|D\psi|\le 1/(t-s)$. Similarly, for $T_1<T_2<T_3$ we say that $\eta$ is a cutoff function for $(T_1,T_3),(T_2,T_3)$ if $\eta$ is a $C^1$ function satisfying $\eta(t)\equiv0$ for $t\le T_1$ and $\eta(t)\equiv1$ if $t\ge T_3$ and $|\eta_t|\le 1/(T_2-T_1)$.

We begin with the following energy estimate for $Du$.
\blemm{dleANSenergy}  We assume that $A,\hat{f}$ satisfy A), F). 
Suppose that $u$ is strong solution  of \mref{e1} on $\Og\times(0,T_0)$. 
Consider any given triple $t_0,T,T'$ satisfying $0<t_0<T<T'\le T_0$ and $p\in [1,n_*/2)$, see the definition \mref{nSGdef} of $n_*$. 

There is a constant $C$, which depends only on the parameters in A) and F), such that for any two concentric balls $B_s,B_t$ with center $x_0\in\bar{\Og}$ and $s<t$ \beqno{keydupANSenergy}\barr{ll}\lefteqn{\sup_{t\in(T,T')}\iidx{\Og_s}{\llg^{-1}(u)|Du|^{2p}}+
	\itQ{Q_{s,t_0}}{|Du|^{2p-2}|D^2u|^2\eta}\le}\hspace{1cm}&\\& C\mathbf{\LLg}^2\itQ{Q_{t,t_0}}{|Du|^{2p+2}\eta}+C((t-s)^{-2}+t_0^{-1})\itQ{Q_{t,t_0}}{|Du|^{2p}}.\earr\eeq
Here, $Q_{t,t_0}=\Og_{t}\times (T-t_0,T')$ and $\eta$ is a cutoff function for $(T-t_0,T'),(T,T')$. \elemm

\bproof  The estimate \mref{keydupANSenergy} results from the energy estimate for $Du$ in \cite[Lemma 3.2]{dleANS} with $W=U=u$ and $\bg(u)=\llg^{-1}(u)$. We differentiated the system in $x$ to obtain
\beqno{ga2zzz} (Du)_t=\Div(A(u)D^2v+A_u(u)DuDu)+D\hat{f}(u,Du).\eeq 
We can
test the above with $\llg^{-1}(u)|Du|^{2p-2}Du\psi^2(x)\eta(t)$ where $\psi$ is a cutoff function for $B_s,B_{t}$. Because $2p<n^*$ and the definition \mref{nSGdef} of $n_*$, it is clear that 
$(2p-2)/(2p) < C_*^{-1}$ so that the spectral gap condition SG) in \cite{dleANS} and \cite[Lemma 3.2]{dleANS} apply here. The only difference is that we keep the function $\eta$ in the integrals of $|Du|^{2p-2}|D^2u|^2$ and $|Du|^{2p+2}$ here and  obtain $$\barr{ll}\lefteqn{\sup_{t\in(T,T')}\iidx{\Og_s}{\llg^{-1}(u)|Du|^{2p}}+
	\itQ{Q_{s,t_0}}{|Du|^{2p-2}|D^2u|^2\eta}\le}\hspace{.5cm}&\\& C\mathbf{\LLg}^2\itQ{Q_{t,t_0}}{|Du|^{2p+2}\eta}+C((t-s)^{-2}+t_0^{-1})\itQ{Q_{t,t_0}}{(1+\llg^{-1}(u))|Du|^{2p}}.\earr$$ Since $\llg(u)$ is bounded from below, the above implies \mref{keydupANSenergy}.
\eproof

Next, we have the following technical result.

\blemm{dleANSprop}  Assume as in \reflemm{dleANSenergy} and that the quantity
\beqno{keydupANSp2}
C_{t_0,T,T'}:=\itQ{\Og\times(T-t_0,T')}{|Du|^{2}} \mbox{ is finite}.\eeq

There is $\mu_0>0$ sufficiently small, in terms of the constants in A) and F), such that if for some positive $R_{\mu_0}$, which may depends on $t_0,T,T'$, such that   \beqno{Keymu01} \mathbf{\LLg}^2\sup_{x_0\in\bar{\Og},t\in (T-t_0,T')}\|u(\cdot,t)\|_{BMO(\Og_{R}(x_0)}^2 \le \mu_0,\eeq
then there are $p>n/2$, an integer $k_0$ and a constant $C^*$ depending only on the parameters of A) and F), $C_{t_0,T,T'}$, $R_{\mu_0}$, and $t_0, T,T'$  such that 
\beqno{keydupANS}\sup_{t\in(T,T')}\iidx{\Og_R}{\llg^{-1}(u)|Du|^{2p}}\le  C^*\mbox{ for any  $R<2^{-k_0}R_{\mu_0}$}.\eeq \elemm

\bproof We follow the argument in the proof of \cite[Proposition 3.1]{dleANS} with $W=U=u$. Suppose that the energy estimate \mref{keydupANSenergy} in \reflemm{dleANSenergy} holds for some $p\ge1$.  We write it as
\beqno{keyiteration1} \ccA(s)+\ccH(s)\le C\mathbf{\LLg}^2\ccB(t)+C[(t-s)^2+t_0^{-1}]\ccC(t), \; 0<s<t,\eeq where the functions $\ccA,\ccH,\ccB$ and $\ccC$ are defined by
$$\ccA(s):=\sup_{t\in(T,T')}\iidx{\Og_s}{\llg^{-1}(u)|Du|^{2p}},\;\ccH(s):=
\itQ{Q_{s,t_0}}{|Du|^{2p-2}|D^2u|^2\eta},$$
$$\ccB(s):= \itQ{Q_{s,t_0}}{|Du|^{2p+2}\eta},\;\ccC(s):=\itQ{Q_{s,t_0}}{|Du|^{2p}}.$$

On the other hand, we apply \reflemm{GNlocalspec} to estimate $\ccB(t)$, the integral of $|Du|^{2p+2}$, on the right hand side of \mref{keyiteration1}. Namely, we let $U=u$ and multiply \mref{GNlocalestspec} by $\mathbf{\LLg}^2\eta$ and integrate the result over $(T-t_0,T')$ to get (recalling the definition of $\mu_0$ in \mref{Keymu01}) 
\beqno{keyiteration2} \mathbf{\LLg}^2\ccB(s) \le \eg\mathbf{\LLg}^2\ccB(t)+C(\eg)\mu_0\ccH(t)+C(\eg)\mu_0(t-s)^{-2}\ccC(t)\quad 0<s<t\le R_{\mu_0}.\eeq

Let us denote $F(t):=\ccB(t)$, $G(t):=\ccH(t)+\mathbf{\LLg}^2\ccB(t)$, $g(t):=\ccC(t)$ and $h(t):=t_0^{-1}\ccC(t)$. The above yields \beqno{Fiter} F(s)\le \eg_0[F(t)+G(t)] +C(t-s)^{-2}g(t)+ Ch(t),\eeq where $\eg_0=\mathbf{\LLg}^2\eg+C(\eg)\mu_0$.
Adding \mref{keyiteration1} and \mref{keyiteration2}, we also have \beqno{Giter}G(s)\le C[F(t)+(t-s)^{-2}g(t)+h(t)].\eeq

It is clear that $\eg_0$ can be small if $\eg$ and then $\mu_0$ are sufficiently small. Thus, if $\mu_0$ is sufficiently small in terms of the constant $C$ of \mref{Giter}, which depends on the constants in A),F), then we can apply a simple iteration argument \cite[Lemma 3.11]{dleANS} to obtain for $0<s<t\le R_{\mu_0}$ $$ F(s)+G(s) \le C[(t-s)^{-2}g(t)+h(t)].$$ Thus, for any $R<R_{{\mu_0}}/2$ we take $t=2R$ and $s=t/2$ and obtain 
\beqno{keydupANSpp}
\itQ{Q_{R,t_0}}{(|Du|^{2p-2}|D^2u|^2+|Du|^{2p+2})}\le C_1(R^{-2}+t_0^{-1})\itQ{Q_{2R,t_0}}{|Du|^{2p}}.\eeq

The above argument shows that if there are $p\ge 1$ and a constant $C(R,t_0)$ such that the energy estimate \mref{keydupANSenergy} holds for $p$ and
\beqno{piterate}\itQ{Q_{2R,t_0}}{\llg^{-1}(u)|Du|^{2p}}\le C(R,t_0),\eeq then this estimate also holds for $p$ being replaced by  any $q\in (p,p+1]$, via \mref{keydupANSpp} and H\"older's inequality. Since \mref{piterate} holds for $p=1$, see \mref{keydupANSp2}, it is now clear that we can repeat the argument $k_0$ times to find  a number $p>n/2$, as long as $2p<n^*$ (so that \mref{keydupANSenergy} holds by \reflemm{dleANSenergy}). We then see that \mref{keydupANSpp} and \mref{piterate} hold for such $p$ and the estimate \mref{keydupANS} follows from the energy estimate for \mref{keydupANSenergy}, with $t=2R,s=R$.
The lemma is proved. \eproof

 The above lemma  made use of a cutoff function $\eta$ for the interval $[T-t_0,T]$ and $[T,T']$ to avoid the dependence on the initial data at $t=0$. This type of result is useful when one wants to discuss the long time dynamics and global attractors of the system.
 
 In order to establish the local and global existence results, we have to provide bounds for $u$ in $\Og\times[0,T_0)$ and allow $t_0=0$. The next lemma considers this case.

\blemm{dleANSpropt0}  Assume as in \reflemm{dleANSprop} with $T=t_0=0$. In addition, we assume that
\beqno{upconti} u \in C([0,T'), L^{2p}(\Og)),\eeq
\beqno{Dupsup} \sup_{t\in[0,T')}\|Du(\cdot,t)\|_{L^{2p}(\Og)}<\infty.\eeq

Then, for the same constant $C^*$ the conclusion \mref{keydupANS} now  reads 
\beqno{keydupANSt0}\sup_{t\in(0,T')}\iidx{\Og_R}{\llg^{-1}(u)|Du|^{2p}}\le   C^*+C\|Du(\cdot,0)\|_{L^{2p}(\Og)}^{2p}\mbox{ for any  $R<2^{-k_0}R_{\mu_0}$}.\eeq
 \elemm

\bproof Thanks to the assumption \mref{Dupsup}, we can let $t_0\to0$ in \mref{keydupANSpp} to see that it also holds for $t_0=0$ and we can also let $T=0$. That is if $\|Du\|_{L^{2p}(Q_{2R,0})}$ is finite then for $Q_{R,0}=\Og_R\times(0,T')$
\beqno{D2uDup} \itQ{Q_{R,0}}{(|Du|^{2p-2}|D^2u|^2+|Du|^{2p+2})} < \infty.\eeq

Using the difference quotience operator $\dg_h$ instead of $D$ in \mref{ga2zzz} in the proof of \reflemm{dleANSenergy}, we obtain
\beqno{ga2zzz1} (\dg_hu)_t=\Div(A(u)D(\dg_h u)+\dg_h(A(u))Du)+\dg_h\hat{f}(u,Du).\eeq 

We test this with $\llg^{-1}(u)|\dg_hu|^{2p-2}\dg_hu\psi^2(x)$ where $\psi$ is a cutoff function for $B_s,B_t$.
Since $u\in C([0,T'),L^{2p}(\Og))$, the energy estimate in \reflemm{dleANSenergy} holds for $t_0=0$ and the operator $D$ is replaced by $\dg_h$. We have
$$\barr{ll}\lefteqn{\sup_{t\in(0,T')}\iidx{\Og_s}{\llg^{-1}(u)|\dg_hu|^{2p}}+
	\itQ{Q_{s,0}}{|\dg_hu|^{2p-2}|D\dg_hu|^2}\le}\hspace{1cm}&\\& C\mathbf{\LLg}^2\itQ{Q_{t,0}}{|Du|^2|\dg_hu|^{2p}}+(t-s)^{-2}\itQ{Q_{t,0}}{|\dg_hu|^{2p}]}+C\iidx{\Og_t}{|\dg_hu(x,0)|^{2p}}.\earr$$

As we now see that the integral in \mref{D2uDup} is finite so that we can let $h$ tend to 0 and obtain a similar energy estimate \mref{keydupANSenergy} for $Du$ with $t_0=0$ and $\eta\equiv1$. Namely,
\beqno{keydupANSenergyt0}\barr{ll}\lefteqn{\sup_{t\in(0,T')}\iidx{\Og_s}{\llg^{-1}(u)|Du|^{2p}}+
	\itQ{Q_{s,0}}{|Du|^{2p-2}|D^2u|^2}\le}\hspace{.5cm}&\\& C\mathbf{\LLg}^2\itQ{Q_{t,0}}{|Du|^{2p+2}}+(t-s)^{-2}\itQ{Q_{t,0}}{|Du|^{2p}]}+C\iidx{\Og_t}{|Du(x,0)|^{2p}}.\earr\eeq

Again, we can argue as in \reflemm{dleANSprop} to treat $\ccB(t)$, the integral of $|Du|^{2p+2}$, on the right hand side and redefine $h(t):=\|Du(\cdot,0)\|_{L^{2p}(\Og)}^{2p}$. The same argument then yields a version of \mref{keydupANSpp} with $t_0=0$. We obtain in particular
\beqno{keydupANSppt0}
\itQ{Q_{R,0}}{|Du|^{2p+2}}\le C_1R^{-2}\itQ{Q_{2R,0}}{|Du|^{2p}}+C_1 \|Du(\cdot,0)\|_{L^{2p}(\Og)}^{2p}.\eeq

The iteration argument after \mref{keydupANSpp} in the proof of \reflemm{dleANSprop} on the power $p$ then gives \mref{keydupANSt0}. This completes the proof.
\eproof

\brem{energyrem}   If the energy estimate \mref{keydupANSenergy} holds with the last term is replaced by the integral of $|Du|^{2p}+\ccF$ for some $\ccF\in L^1_{loc}(Q)$ then the argument in \reflemm{dleANSprop} yields \mref{keydupANSpp} with the integral of $|Du|^{2p}$ on right hand side being  replaced by that of $|Du|^{2p}+\ccF$ over $Q_{2R,t_0}$. Similar observation applies to \reflemm{dleANSpropt0}.

\erem

We are now ready to provide the proof of the main theorem. 

{\bf Proof of \reftheo{dleNL-mainthm}:} By \reflemm{Tmaplem}, the map $T_\sg: \mX\to\mX$ defined by \mref{Tmapdef} is compact. In order to apply the Leray-Schauder theorem and show that there is a fixed point $u$ for $\sg=1$, which the solution of \mref{e1}, we need only to provide a uniform bound for the fixed points of $T_\sg$. To this end, for any $\sg\in[0,1]$ we consider a fixed point $u$ of $T_\sg(u)=u$. 

Since $u\in\mX$, $u$ is a bounded weak solution and continuous so that \cite[Theorems 2.1 and 3.1]{GiaS} apply and yield that $Du$ is locally bounded in $Q=\Og\times(0,T_0)$. It is then well known that $D^2u$ exists and $D^2u\in L_{loc}(Q)$ and thus $u$ is a strong solution in $Q$. 

We will apply \reflemm{dleANSpropt0} here. First of all, the continuity assumption \mref{upconti} of the lemma is clear because $u\in \mX$. 
Next, for any $q=2p\in(n,p_0)$ we show that  $\|u(\cdot,t)\|_{W^{1,q}(\Og)}$ is bounded $[0,T_0)$ to verify \mref{Dupsup}. For any $h>0$ and any function $w$ we denote by $w_{(h)}=\fg_h*w$ the mollifier/regularizer of $w$. We have for any $f\in L^{q'}(\Og)$ $$\barr{lll}\iidx{\Og}{D(u(x,t)_{(h)})f(x)} &=&\iidx{\Og}{(Du(x,t))_{(h)}f(x)}=\iidx{\Og}{Du(x,t)f_{(h)}(x)}\\&=&
\iidx{\Og}{u(x,t)Df_{(h)}(x)}\to \iidx{\Og}{u(x,0)Df_{(h)}(x)},\earr$$ as $t\to0$ because $u\in\mX$. The last term in the above is bounded by $\|u(\cdot,0)\|_{W^{1,p}(\Og)}\|f_h\|_{L^{q'}(\Og)}$ or $C\|U_0\|_{W^{1,p_0}(\Og)}\|f\|_{L^{q'}(\Og)}$. By the Uniform Boundedness Principle, noting that $Du(\cdot,t)\in L^q(\Og)$ for each $t>0$, we see that $\|Du_{(h)}(\cdot,t)\|_{L^q(\Og)}$ is uniformly bounded for all $h>0$ and $t\in[0,T_0)$. Let $h\to0$ we derive that $\|Du(\cdot,t)\|_{L^{2p}(\Og)}$ is bounded for $t\in[0,T_0)$. Thus, for each fixed point $u$ of $T_\sg$ the condition \mref{Dupsup} holds.

Hence, from the assumptions M.1) and M.2), \reflemm{dleANSpropt0} can apply here to provide uniform constants $C^*,R_1$ depending only on the parameters of A) and F), $C_{T_0}$, $R_{\mu_0}$, and $\|DU_0\|_{L^{p_0}(\Og)}$  such that if $p<p^*=\frac12\min\{n^*,p_0\}$ then
\beqno{keydupANSt0zz}\sup_{t\in(0,T_0)}\iidx{\Og_{R_1}}{\llg^{-1}(u)|Du|^{2p}}\le   C^*.\eeq 

 From the definition of $r^*$ it is clear that we can choose $p,p_1$ such that $n<p_1<p<p^*$ and $r^*=p_1/(p-p_1)$.  As $r^*=p_1/p(p/p_1)'$, by H\"older's inequality we have
\beqno{keydupANSt0zzzz}\iidx{\Og_{R_1}}{|Du|^{2p_1}}\le\|\llg(u)\|_{L^{r^*}(\Og_{R_1})}\left( \iidx{\Og_{R_1}}{\llg^{-1}(u)|Du|^{2p}}\right)^{p_1/p}.\eeq 
From the assumption L) on $\llg(u)$ and   \mref{keydupANSt0zz},  the right hand side of \mref{keydupANSt0zzzz} will be bounded uniformly for all $\sg\in[0,1]$.

We then have a uniform bounded for $\|u\|_{W^{1,q}(\Og)}$. As $q=2p_1>n$, by Sobolev's embedding theorem, we see that $\|u\|_\mX\le M$ for some constant $M$ and all $\sg\in[0,1]$.  The Leray Schauder theory then applies to provide a fixed point $u=T_1(u)$. By  \reflemm{Uniquelem} at the end of this paper, which deals with uniqueness of continuous weak solutions, this fixed point is the unique strong solution of the system \mref{e1}. The proof is then complete.
\eproof

{\bf Proof of \refcoro{dleNL-maincoro}:} For any positives $T_1,T_2$ with $T_1<T_2$, since we are assuming that the conditions of \reftheo{dleNL-mainthm} hold for all $T_0>0$, we see that \mref{e1} has strong solutions $u_i$'s on $Q_i=\Og\times(0,T_i)$, $i=1,2$. We also showed that there is a constant $C_i$ such that $\sup_{[0,T_i)}\|u_i\|_{W^{1,2p}(\Og)}\le C_i$ for some $p>n/2$. By  \reflemm{Uniquelem} on the uniqueness of the weak solutions satisfying this boundedness, we see that $u_1\equiv u_2$ on $Q_1$. Thus, \mref{e1} has a unique strong solution which exists globally. \eproof

{\bf Proof of  \refcoro{dlec1coro}:} We just need to show that the assumption D) implies M.1) and M.2). It is clear that D) yields M.2). To verify M.1)
we argue by contradiction.  If this is not the case then there are sequences $\{x_n\}\subset\bar{\Og}$, $\{\sg_n\}\subset [0,1]$, $\{t_n\}\subset (0,T_0)$, $\{r_n\}$, $r_n\to0$ and a sequence of strong solutions $\{u_{\sg_n}\}$ such that for $U_n(\cdot)=u_{\sg_n}(\cdot,t_n)$ $$\|U_{n}\|_{BMO(B_{r_n}(x_n)\cap\Og)}>\eg_0 \mbox{ for some $\eg_0>0$}.$$ 

By D) we see that the sequence $\{U_n\}$ is bounded in $W^{1,n}(\Og)$. We can then assume that $U_n$ converges weakly to some $U$ in $W^{1,2}(\Og)$ and strongly in $L^2(\Og)$.  We then have $\|U_n\|_{BMO(B_R\cap\Og)}\to \|U\|_{BMO(B_R\cap\Og)}$ for any given ball $B_R$. It is easy to see that $U\in W^{1,n}(\Og)$ and by Poincar\'e's inequality $U$ is VMO and $\|U\|_{BMO(B_R\cap\Og)}<\eg_0/2$ if $R$ is sufficiently small. The number $R$ is independent of $\llg_0\ge1$ because $\|U\|_{W^{1,n}(\Og)}$ is independent of $\llg_0$. Furthermore, we can assume also that $x_n$ converges to some $x\in\bar{\Og}$. Thus, for large $n$, we have $r_n<R/2$ and $x_n\in B_{R/2}(x)$. Then, for large $n$,  $B_{r_n}(x_n)\subset B_{R}(x)$ and $$\|U_n\|_{BMO(B_{r_n}(x_n)\cap\Og)} \le \|U_n\|_{BMO(B_R(x)\cap\Og)}\le \|U\|_{BMO(B_R(x)\cap\Og)}+\eg_0/2<\eg_0.$$ We obtain a contradiction. Thus, M.1) holds and the proof is complete. \eproof

{\bf Proof of  \refcoro{dlec2coro}:} We need only show that D) and L') together imply L).
Let $u$ be any strong solution of \mref{e1famzzz} and $\llg(u)$ satisfy \mref{llgullgeg}. There are positives $s_0,C_0$ and such that \beqno{llgs0}\|\llg^{s_0}(u)\|_{L^{1}(\Og)}\le C_0(T_0). \eeq 
We will show that for any $r>1$ there is a constant $C$ depending on $C_0,s_0,r,|\Og|,T_0$ and $\|u\|_{W^{1,n}(\Og)}$ such that \beqno{llgs0r}  \|\llg(u)\|_{L^{r}(\Og)}\le C.\eeq

 We choose and fix $s>0$ and $p\in(1,n)$ such that $sp_*=s_0$, where $p_*=np/(n-p)$. Then \mref{llgs0} implies \beqno{llgs} \|\llg^{s}(u)\|_{L^{p_*}(\Og)}\le C_0^\frac{1}{p_*}(T_0).\eeq 

Define $g(\cdot)=\llg^{s+\eg_0}(u(\cdot,t))$. The definition of $\BLLg_1$ in \mref{llgullgeg} gives
$$|Dg|\le C(s)\frac{|\llg_u|}{\llg^{1-\eg_0}(u)}\llg^s(u)|Du| \le C(s)\BLLg_1 \llg^s(u)|Du|.$$ Hence,  by H\"older's inequality, $\|Dg\|_{L^p(\Og)} \le C\|\llg^s(u)\|_{L^{p_*}(\Og)}\|Du\|_{L^n(\Og)}$. This and \mref{llgs} and D) imply some $C(T_0)$ such that $\|Dg\|_{L^{p}(\Og)}\le C(T_0)$. Using H\"older's inequality, we have $\|g\|_{L^1(\Og)}\le C\|\llg^s(u)\|_{L^{p_*}(\Og)}^{(s+\eg_0)/s}\le CC_0^{1+\eg_0/(p_*s)}(T_0)$. Hence, $$\|g\|_{W^{1,p}(\Og)}\le C(T_0)+CC_0^{1+\eg_0/(p_*s)}(T_0).$$
By Sobolev's embedding theorem, $\|g\|_{L^{p_*}(\Og)}$ is bounded. From the definition of $g$, we can find a constant $\bar{C}(T_0)$ such that $\|\llg^{s+\eg_0}(u)\|_{L^{p_*}(\Og)}\le \bar{C}(T_0)$. Thus, there is a constant $C_1(T_0)$ such that
$$\|\llg^{s_0+p_*\eg_0}(u)\|_{L^{1}(\Og)} \le C_1(T_0).$$

This shows that if \mref{llgs0} holds for some $s_0$ then it also holds for $s_0$ being $s_0+p_*\eg_0$ and new constant $C_1(T_0)$.  It is then clear that we can repeat this argument  to see that $\|\llg^{s_0+kp_*\eg_0}(u)\|_{L^{1}(\Og)} \le C_k(T_0)$ for all integers $k$ and some $C_k(T_0)$. This fact and a simple use of H\"older's inequality show that L) holds. The proof is complete. \eproof

\brem{uDuagreminx} By \mref{keydupANSt0zzzz}, $u$ is H\"older in $x$. We can show that $u$ is also H\"older continuous in $x,t$. Indeed, \mref{keydupANSpp} with $p=1$ shows that $|D^2u|^2, |Du|^4$ are in $L^1(Q)$. We obtain from the system of $u$ and a simple use of H\"older's inequality that
$$\|u_t\|_{L^1(Q)}\le \|A(u)\|_{L^2(Q)}\|D^2u\|_{L^2(Q)}+\|A_u(u)\|_{L^2(Q)}\||Du|^2\|_{L^2(Q)}+\|\hat{f}\|_{L^1(Q)}.$$

Since $|A(u)|,|A_u(u)|\le \llg(u)$ and $\hat{f}$ has linear growth in $Du$, the right hand side is finite and bounded by a constant independent of $\llg_0$ (using \mref{keydupANSpp} with $p=1$ and then \mref{keydupANSt0zzzz} to see that $\|D^2u\|_{L^2(Q)}\le C\llg_0^{-1}$). Thus, $u_t$ belongs to $L^1(Q)$. It is is well known that if $u$ is H\"older continuous in $x$ and $u_t$ is in $L^1(Q)$ then $u$ is H\"older in $x,t$ (see \cite[Lemma 4]{NSver}).\erem

\section{The general reaction term and the proof of \reftheo{Genthm}}\label{genreaction}\eqnoset

We provide the proof of \reftheo{Genthm} in this section. Recall the assumption \mref{DFDu0} of \reftheo{Genthm} \beqno{DFDu} |\partial_\zeta \hat{f}(u,\zeta)|\le C\llg^\frac12(u),\; |\partial_u \hat{f}(u,0)|\le C\llg(u).\eeq

We then define $F_\sg(u,\zeta):=\hat{f}(\sg u,\sg\zeta)$ and \beqno{Bfdef}B_\sg(u,\zeta):=\int_{0}^{1}\partial_\zeta F_\sg(u,t\zeta)\,dt,\quad f_\sg(u):=\int_{0}^{1}\partial_u F_\sg(tu,0)\,dt.\eeq

Fixing $h>0$, for any given $u\in\mX$ and $\sg\in[0,1]$ we consider the following linear system
\beqno{Tmapdefa}\left\{\barr{l} w_t=\Div(A(\sg u)Dw)+B_\sg^{(h)}(u)Dw+f_\sg(u)w+F_\sg(0,0)\quad (x,t)\in Q, \\w(x,0)=U_0(x)\quad  x\in\Og\\\mbox{$w=0$ or $\frac{\partial w}{\partial \nu}=0$ on $\partial \Og\times(0,T_0)$}, \earr\right.\eeq
where $B_\sg^{(h)}(u):=B_\sg(u,D(u_{(h)}))$ with $u_{(h)}$ being the mollifier/regularizer of $u$.

By \mref{DFDu}, $B_\sg^{(h)}(u)$ is bounded if $u\in\mX$. Therefore, we can follow \reflemm{Tmaplem} and use \mref{Tmapdefa} to define the compact map $T_\sg:\mX\to\mX$ by $T_\ag(u)=w$, the weak solution of \mref{Tmapdefa}.

We now consider a fixed point $u$ of $T_\sg$ and provide a uniform estimate for $\|u\|_\mX$. Note that the results in \cite{GiaS} still applies to the system so that $u$ is a strong solution. We will apply \reflemm{dleANSpropt0} again. The first task is to establish the energy estimate \mref{keydupANSenergy} in \reflemm{dleANSenergy} for $Du$ with the new $\hat{f}(u,Du)$ is now
$$F(u,Du)=B_\sg^{(h)}(u)Du+f_\sg(u)u+F_\sg(0,0).$$

 Again, we differentiate the system with respect to $x$ and test the result with $Du|^{2p-2}Du\psi^2\eta$. We go back to the proof of \cite[Lemma 3.2]{dleANS} and \cite[Remark 3.4]{dleANS} to see that we have to deal with the integral of $D(B_\sg^{(h)}(u)Du)|Du|^{2p-2}Du\psi^2\eta$. Integrating by parts, we have
$$\itQ{Q}{D(B_\sg^{(h)}(u)Du)|Du|^{2p-2}Du\psi^2\eta}=-\itQ{Q}{B_\sg^{(h)}(u)DuD(|Du|^{2p-2}Du\psi^2)\eta}.$$

Expanding $D(|Du|^{2p-2}Du\psi^2)$ and using \mref{DFDu} we can see easily that the energy estimate for $Du$ still holds. Therefore, \reflemm{dleANSpropt0} applies to give a uniform estimate for $\|u\|_{W^{1,p}(\Og)}$ for some $p>n$.  The same fixed point argument now gives a fixed point $U^{(h)}$ of $T_1$, solving
\beqno{Uheqn}(U^{(h)})_t=\Div(A(U^{(h)})DU^{(h)})+B_1(U^{(h)},DU_{(h)}^{(h)})DU^{(h)}+f_1(U^{(h)})U^{(h)}+F_1(0,0).\eeq

Thus, for each $h>0$ we obtain a strong solution $U^{(h)}$.
These solutions are uniformly bounded in $L^\infty([0,T_0),W^{1,p}(\Og))$, with $p>n$, so that there exists a sequence $h_n\to0$ and $U_n=U^{(h_n)}$ converges in $\mX$ to some $u$. 

Moreover, for any $t_0>0$ and any compact set $K\subset \Og\times(t_0,T_0)$ we can use the results in \cite{GiaS} to see that $\{DU_n\}$ is uniformly bounded in $C^{\ag,\ag/2}(K)$. Moreover, $D^2U_n$ is bounded uniformly in $L^2(K)$, see \mref{keydupANSpp} with $p=1$. Thus, there is a subsequence $\{U_{n_k}\}$ such that $DU_{n_k}$ converges  to $Du$ in $C(K)$. Letting $h_n\to0$ in \mref{Uheqn}, we see that $u$ is a continuous weak solution to the following system in   
$\Og\times(t_0,T_0)$
$$u_t = \Div(A(u)Du) + B_1(u,Du)Du + f_1(u)u+F_1(0,0).$$ On the other hand, from the definition \mref{Bfdef}, $F_\sg(u,\zeta)=B_\sg(u,\zeta)\zeta+f_\sg(u)u+ F_\sg(0,0)$ and because $\hat{f}(u,Du)=F_1(u,Du)$, we see that $u$ solves \beqno{genFsys}u_t = \Div(A(u)Du) + \hat{f}(u,Du).\eeq

 Since $u,Du$ are bounded and continuous, $u$ is a desired strong solution in $\Og\times(0,T_0)$. To complete the proof of the theorem, we need only show that $u$ is unique. Note that the uniform assumptions in the theorem implies that any strong solution $u$ is also a weak solution satisfying for some constants $C^*$ and $p>n/2$
\beqno{unikeydupANSt0zz}\sup_{t\in(0,T_0)}\iidx{\Og}{|Du|^{2p}}\le   C^*.\eeq

The uniqueness of $u$ then follows from \reflemm{Uniquelem} below.

We consider the following general quasilinear system.
\beqno{gensys}\left\{\barr{l} u_t=\Div(A(u,Du))+F(u,Du)\quad (x,t)\in Q=\Og\times(0,T_0),\\u(x,0)=U_0(x)\quad  x\in\Og\\\mbox{$u=0$ or $\frac{\partial u}{\partial \nu}=0$ on $\partial \Og\times(0,T_0)$}. \earr\right.\eeq

Here, $A,F$ are maps from $\RR^m\times\RR^{mn}$ into $\RR^m$.
We assume the following conditions.
\bdes\item[A.1)] (Monotonicity condition) There is $\llg_0>0$ such that for any $u\in\RR^m$ and $\zeta_1,\zeta
_2\in\RR^{nm}$
\beqno{AAcond1}\llg_0|\zeta_1-\zeta_2|^2 \le \myprod{A(u,\zeta_1)-A(u,\zeta_2),\zeta_1-\zeta_2}.\eeq

\item[A.2)] (Lipschitz condition) There is a continuous function $C$ on $\RR^m\times\RR^m$ such that for any $u_1,u_2\in\RR^m$ and $\zeta_1,\zeta
_2\in\RR^{nm}$
\beqno{AAcond2} |A(u_1,\zeta_1)-A(u_2,\zeta_1)|\le C(u_1,u_2)|u_1-u_2||\zeta_1|.\eeq
\beqno{FFcond2}|F(u_1,\zeta_1)-F(u_2,\zeta_2)| \le C(u_1,u_2)(|u_1-u_2||\zeta_2|+|\zeta_1-\zeta_2|).\eeq

\edes

It is easy to see that \mref{AAcond1} and \mref{AAcond2} are satisfied if the matrix $A_\zeta$ is positive definite and $A$ is Lipschitz in $u$ and has a linear growth in $\zeta$. Concerning \mref{FFcond2}, we note that $$|F(u_1,\zeta_1)-F(u_2,\zeta_2)|\le|F(u_1,\zeta_2)-F(u_2,\zeta_2)|+|F(u_1,\zeta_1)-F(u_1,\zeta_2)|$$ so that \mref{FFcond2} can be verified if $F(u,Du)$ is Lipschitz in $u$ and $Du$.

\blemm{Uniquelem} Assume A.1) and A.2). Suppose that there are constants $C^*$ and $p>n/2$ such that {\em any} weak solution $u$ of \mref{gensys} on $Q=\Og\times(0,T_0)$ satisfies
\beqno{uniboundforDu}\sup_{t\in(0,T_0)}\iidx{\Og}{|Du|^{2p}}\le   C^*.\eeq

Then \mref{gensys} has at most one weak solution on $Q$.
\elemm

\bproof
Let $u_1,u_2$ be two weak solutions to \mref{gensys} and $w=u_1-u_2$. Substracting the two systems and rearranging, we have
$$w_t = \Div(I_1 +I_2) + I_3,$$ where $I_1:=A(u_1,Du_1)-A(u_1,Du_2)$,
$$I_2:=A(u_1,Du_2)-A(u_2,Du_2)\mbox{ and }I_3:=F(u_1,Du_1)-F(u_2,Du_2).$$

Multiplying the system for $w$ with $w$, integrating over $Q_t=\Og\times(0,t)$ and integrating by parts in $x$, since $w=0$ at $t=0$, we obtain
$$\iidx{\Og}{|w(\cdot,t)|^2} +\itQ{Q_t}{\myprod{I_1,Dw}}= \itQ{Q_t}{[-\myprod{I_2,Dw}+\myprod{I_3,w}]},$$ 

In the sequel, we will use the assumptions A.1) and A.2) with $\zeta_1=Du_1$, $\zeta_2=Du_2$. In particular, using A.1) in $\myprod{I_1,Dw}$, we get
\beqno{wtest}\iidx{\Og}{|w(\cdot,t)|^2}+\llg_0\itQ{Q_t}{|Dw|^2} \le \int_0^t[J_2(s)+J_3(s)]\,ds,\eeq where
\beqno{J2}J_2(s):=\iidx{\Og}{|A(u_1,Du_2)-A(u_2,Du_2)||Dw|},\eeq
\beqno{J3}J_3(s):=\iidx{\Og}{|F(u_1,Du_1)-F(u_2,Du_2)||w|}.\eeq

The assumption \mref{uniboundforDu} implies a constant $C( C^*)$ such that \beqno{UNIDu}\sup_{s\in(0,T_0)}\|u_i(\cdot,s)\|_{L^{\infty}(\Og)},\;\sup_{s\in(0,T_0)}\|Du_i(\cdot,s)\|_{L^{2p}(\Og)}\le  C( C^*)\quad i=1,2.\eeq

Concerning $J_2$, we use \mref{AAcond2} and the boundedness of $u_1,u_2$ in the above and H\"older's inequality to get for $s\in(0,T_0)$
$$J_2(s)\le C(C^*)\||w||Du_2|\|_{L^2(\Og)}\|Dw\|_{L^2(\Og)}.$$By H\"older's inequality and \mref{UNIDu} again, we have
\beqno{wDu2}C\||w||Du_2|\|_{L^2(\Og)}\le C\|w\|_{L^{2p'}(\Og)}\||Du_2|\|_{L^{2p}(\Og)} \le  C(C^*)\|w\|_{L^{2p'}(\Og)},\eeq where $2p'=2p/(p-1)$. Because $p>n/2$, $2p'<2n/(n-2)$ so that we can use the Poincar\'e-Sobolev inequality to have for any given $\eg>0$
\beqno{wDu2z}C(C^*)\|w\|_{L^{2p'}(\Og)} \le \eg\|Dw\|_{L^2(\Og)}+C(\eg,C^*)\|w\|_{L^2(\Og)}.\eeq

We now see that $J_2(s)\le C(C^*)\|w\|_{L^{2p'}(\Og)}\|Dw\|_{L^{2}(\Og)}$. Using \mref{wDu2z} and Young's inequality, we have
\beqno{J2est}J_2(s)\le\eg\|Dw\|_{L^2(\Og)}^2+C(\eg,C^*)\|w\|_{L^2(\Og)}^2\quad \mbox{for all $s\in(0,T_0)$}.\eeq

Similarly, from \mref{FFcond2}
$J_3(s)\le C(C^*)[\||w||Du_2|\|_{L^2(\Og)}+\|Dw\|_{L^2(\Og)}]\|w\|_{L^2(\Og)}$ and we can use \mref{wDu2}, \mref{wDu2z} and then Young's inequality to see that
\beqno{J3est}\barr{lll} J_3(s)&\le&C(C^*)[\|w\|_{L^{2p'}(\Og)}+\|Dw\|_{L^2(\Og)}]\|w\|_{L^2(\Og)}\\&\le& \eg\|Dw\|_{L^2(\Og)}^2+C(\eg,C^*)\|w\|_{L^2(\Og)}^2 .\earr\eeq 

Therefore, using \mref{J2est} and \mref{J3est} in \mref{wtest}, we obtain
$$\iidx{\Og}{|w(\cdot,t)|^2}+\llg_0\itQ{Q_t}{|Dw|^2} \le\int_0^t(\eg\|Dw(\cdot,s)\|_{L^2(\Og)}^2+C(\eg,C^*)\|w(\cdot,s)\|_{L^2(\Og)}^2)\,ds.$$
Choosing $\eg$ small, we the obtain the following integral Gronwall inequality. 
$$\|w(\cdot,t)\|_{L^2(\Og)}^2 \le C(C^*)\int_0^t\|w(\cdot,s)\|_{L^2(\Og)}^2\,ds.$$

Since $\|w(\cdot,0)\|_{L^2(\Og)}^2=0$, we obtain $w\equiv0$, i.e., $u_1\equiv u_2$, on $Q$. The proof is complete. \eproof

\bibliographystyle{plain}

\end{document}